\documentclass[12pt]{article}
\usepackage{amssymb,amsmath}
\def\C{{\mathbf{C}}}
\def\bC{{\mathbf{\overline{C}}}}
\def\P{{\mathbf{P}}}
\def\N{{\mathfrak{N}}}
\def\n{{\mathfrak{n}}}
\newtheorem{lemma}{Lemma}
\newtheorem{theorem}{Theorem}
\begin{document}
\author{Alexandre Eremenko\thanks{Supported by NSF grant DMS-1067886}}
\title{Normal holomorphic maps from $\C^*$ to a projective space}
\maketitle
\begin{abstract}
A theorem of A. Ostrowski describing meromorphic functions $f$ such that
the family $\{ f(\lambda z):\lambda\in\C^*\}$ is normal,
is generalized to holomorphic maps from $\C^*$ to a projective space.
\end{abstract}
\maketitle

Let $f:\C^*\to\P^n$ be a holomorphic curve and 
\begin{equation}\label{1}
F=(g_0,g_1,\ldots,g_n)
\end{equation}
some homogeneous representation of $f$. This means that $g_j$ are analytic functions
in $\C^*$
without common zeros. When $n=1$, we have $\P^1=\bC$,
and $f$ can be identified
with a meromorphic function $g_1/g_0$ in $\C^*$.

There is a conformal Riemannian metric with the line element 
\begin{equation}\label{metric}
|dz|/|z|
\end{equation}
on $\C^*$
which is invariant with respect to conformal automorphisms
$z\mapsto\lambda z,\;\lambda\neq 0$. The punctured plane with this metric
is isometric to a cylinder of infinite length and circumference $2\pi$. 

A holomorphic curve in $\C^*$ is called {\em normal}
if it is uniformly continuous
with respect to this metric
(\ref{metric}) and the Fubini-Study metric in $\P^n$.
An equivalent property is that $\{ z\mapsto f(\lambda z):\lambda\in\C^*\}$
is a {\em normal family}: every sequence of these maps has a subsequence
which converges uniformly on compacts with respect to the Fubini--Study metric
to a holomorphic map $\C^*\to\P^n$.

A normal holomorphic curve in $\C^*$ has genus zero, which means that it has 
a {\em canonical} homogeneous representation
\begin{equation}\label{co}
g_j(z)=A_jz^{m_j}\prod_{|z_{j,k}|<1}(1-z_{j,k}/z)\prod_{|z_{j,k}|\geq 1}(1-z/z_{j,k}),
\end{equation}
where $z_{j,k}$ are the zeros of $g_j$, $A_j$ is a constant,
and $m_j$ is an integer. We can and will always assume that $\min_jm_j=0$,
which defines these integers uniquely. It is clear that $z_{j,k}$ 
are uniquely defined by $f$, and $A_j$ are defined up to a common multiple.

In fact, 
$$A_j=\exp\left(\frac{1}{2\pi}\int_{-\pi}^\pi\log|g_j(e^{i\theta})|d\theta
\right).$$

A. Ostrowski \cite{O} considered normal meromorphic functions ($n=1$),
and completely characterized them in terms of parameters of the 
canonical representation,
see also \cite[Ch. VI]{M}, for an exposition of this work.

In this paper, the result of Ostrowski is extended to arbitrary dimension~$n$.

We begin with a reformulation of convergence of curves in terms of homogeneous coordinates.
Consider a sequence $F_k=(g_{k,0},\ldots,g_{k,n})$ of $(n+1)$-tuples of
holomorphic functions in an arbitrary region $D$.
We assume that coordinates of each $F_k$ have no zeros common to all of them.
Then we have a sequence of holomorphic curves $f_k:D\to\P^n$.

\begin{lemma}\label{lemma1}
The sequence $(f_k)$ converges with respect to the Fubini--Study metric,
uniformly on compacts in $D$,
if and only if the following two
conditions are satisfied:

\noindent
(i) For every compact $K\subset D$,
there exist functions $h_k$ holomorphic on $K$, having no zeros,
and such that for every $j$ there exists a limit
uniform on $K$, with respect to the Euclidean metric
in $\C$: 
\begin{equation}\label{a}
\lim_{k\to\infty} h_kg_{k,j}=g_{\infty,j},
\end{equation}
and

\noindent
(ii) If (\ref{a}) is satisfied with some $h_k$ as in (i)
then 
functions $g_{\infty,0},\ldots,g_{\infty,n}$ have no common zeros.
\end{lemma}

Some, but not all, functions $g_{\infty,j}$ may be identically equal to $0$.
Condition (ii) means that the rest of them have no common zeros.
\vspace{.1in}

{\em Proof.}
It is evident that (i) and (ii) imply convergence of $(f_k)$.

In the other direction, let $f_k\to f$ converge.  Let
$F=(g_0,\ldots,g_n)$ be a homogeneous representation of $f$.
Let
$I\subset\{0,\ldots,n\}$ be the set of indices for which $g_j\not\equiv 0$.
There exist $\delta>0$ and a finite open covering $\{ D_j:j\in I\}$ of $K$
such that $|g_j(z)|\geq \delta,\; z\in D_j,\; j\in I$. This implies that
Fubini--Study distance from $f(D_j)$ to the hyperplane $w_j=0$ is
positive, so
$g_{k,j}$ are free from zeros on $D_j$ when $k$ is large enough.

On $D_j$ we define $h_{k,j}=g_j/g_{k,j}$. These functions are
holomorphic and zero-free on $D_j$ when $k$ is large enough.
We have $g_{k,i}/g_{k,j}\to g_i/g_j$ uniformly on $D_i\cap D_j$.
Let $p_{k,i,j}=h_{k,i}/h_{k,j}$ on $D_i\cap D_j$.
Then 
\begin{equation}\label{11}
p_{k,i,j}\to 1
\end{equation}
uniformly on $D_i\cap D_j$ and we have the cocycle condition
\begin{equation}\label{cocycle}
p_{k,i,j}p_{k,j,l}p_{k,l,i}=1
\end{equation}
on triple intersections.
In view of (\ref{11}) we can define the principal branches of $\log p_{k,i,j}$
on the double intersections. Then there exist holomorphic functions
$\phi_{k,j}$ on $D_j$ such that
$\log p_{k,i,j}=\phi_{k,i}-\phi_{k,j}$ on $D_i\cap D_j$,
and $\phi_{k,j}\to 0$
as $k\to\infty$, on $D_j$ (we may need to shrink the $D_j$ a little
at this step). See \cite{Hor},
theorems 1.2.2, $1.4.3^\prime, 1.4.4$ and $4.4.2$.
Now we set
$h_k=h_{k,j}\exp(-\phi_{k,j})$ and these $h_k$ do the job. This proves (i).

The functions $h_k$ we constructed have property (ii).
Now we show that (ii) must hold for every sequence of functions $h_k$
as in (i).
Let $I^\prime=\{ j:g_{\infty,j}\equiv 0\}$,
and $I=\{0,\ldots,n\}\backslash I^\prime.$ Suppose that
$z_0$ is a common zero of $g_{\infty,j},\; j\in I.$
Then for every $\epsilon>0$ there
is a closed disc $G$ centered at $z_0$
such that for $j\in I$ and $k$ large enough,
each function $h_kg_{j,k}$ has a zero in $G$.
Let
$$M_k=\max\{ |h_kg_{k,j}(z)|:z\in G, j\in I\}.$$
Then $M_k$ is bounded from below as $k\to\infty$ by a constant
that depends only on $G$. So we have
$$\max\{ |h_kg_{k,j}(z)|:z\in G, j\in I^\prime\}=o(M_k).$$
This means that for some points $z_k\in G$, $f_k(z_k)$ tends to the subspace
$$H_{I'}=\{(w_0,\ldots,w_n):w_j=0, j\in I^\prime\}.$$
On the other hand,
$f_k(z),\; z\in G$
visits every hyperplane $H_j$ defined by $w_j=0$ for $j\in I$.
As these hyperplanes and subspace $H_{I^\prime}$ have empty intersection,
diameter of $f_k(G)$ must be greater then a positive constant independent
of $k$. This contradicts our assumption that $G$ is a disc of radius
$\epsilon$ which can be arbitrarily small, because the sequence
$(f_k)$ is equicontinuous.
This completes the proof of Lemma~\ref{lemma1}.
\vspace{.1in}

For the future use, we need a restatement of condition (ii) which does not
involve the limit functions $g_{\infty,j}$.

Suppose that (i) holds. Then (ii) is equivalent to the following condition:

{\em There exist $C$ and $\delta$ (depending on $K$ and $(f_k)$) such that
such that for every disc $D(z_0,\delta)$  with $z_0\in K$, and for every
$I\subset\{0,\ldots,n\}$, whenever all $g_{k,j},\; j\in I$ have zeros
in $D(z_0,\delta)$, we have
\begin{equation}\label{cond2}
\max_{0\leq j\leq n}|h_kg_{k,j}(z_0)|\leq \max_{j\in I^\prime}|h_kg_{k,j}(z_0)|+C,
\end{equation}
where $I^\prime=\{0,\ldots,n\}\backslash I$.}
\vspace{.1in}

The equivalence of (ii) and (\ref{cond2}), assuming that (i) holds, has been established
in the proof of the second part of Lemma~\ref{lemma1}.
\vspace{.1in}

For a function $g$ holomorphic in a ring $\{ z:r_1<|z|<r_2\}$, we define
$$N(r,g)=\frac{1}{2\pi}\int_{-\pi}^\pi\log|g(re^{i\theta})|d\theta,$$
and 
$$\N(t,g)=N(e^t,g).$$
It is well-known that this function $\N$ is convex on $(r_1,r_2)$,
and piecewise-linear. It is linear (affine) on an interval $(a,b)$ if $g$ has
no zeros in the ring $\{ z:e^a<|z|<e^b\}$,
and the derivative $\N'$ has a jump $k$
at the point $t$ if $g$ has $k$ zeros on the circle $|z|=e^t$.
All this follows from the Jensen formula.

\begin{lemma}\label{lemma2} Using the notation of Lemma 1, suppose that $F_k$
are defined in a ring $\{ r_1<|z|<r_2\}$ and that $(f_k)$ converges
to a limit, so that (i) and (ii) hold.
Then there exist linear functions $\ell_k(t)=a_kt+b_k$, such
that for every $j\in[0,n]$ the limit
$$\lim_{k\to\infty}\N(t,g_{k,j}-\ell_k)<+\infty$$
exists, possibly identically equal to $-\infty$, uniformly on every
interval $[a,b]$ such that $\log{r_1}<a<b<\log r_2.$
\end{lemma}
\vspace{.1in}

{\em Proof.} Indeed, the functions $h_k$ of Lemma~\ref{lemma1} are zero-free
in the ring, so $\N(t,h_k)$ are linear functions.
\vspace{.1in}

Let us fix some interval $(-a,a)$ and consider $(n+1)$-tuples
$$\Phi=(\phi_0,\phi_1,\ldots,\phi_n)$$ of convex functions on $(-a,a)$.
We say that a sequence $\Phi_k$ of such tuples {\em converges uniformly}
if for each $j$ the coordinates $\phi_{k,j}$ converge uniformly
on compact subintervals to finite convex functions, or to identical $-\infty$.
A family of such $(n+1)$-tuples of convex functions is called {\em normal}
if every sequence contains a subsequence that converges uniformly, 
but not all coordinates converge to $-\infty$.

If $n=0$ and we are dealing with a family of convex functions,
then the criterion
of normality is that the family is uniformly bounded from above on each
compact subinterval.
The limit function is finite if in addition the functions
of the sequence are bounded from below at some point. These
statements are well-known and
easy to prove.

An equivalent criterion of normality with all limit functions finite
for  $n=0$
is that {\em all functions are bounded
at some point, their derivatives are bounded 
at the same point, and the total jump of the derivatives is bounded
on each compact subinterval.}

From this, it is easy to derive a criterion for every $n>0$:
{\em for normality of a family of
$(n+1)$-tuples of convex functions,
it is necessary and sufficient that the functions
$\phi=\max_j\phi_j$ form a normal family with $n=0$.}

\begin{lemma}\label{lemma3} 
Let $X=\{\Phi\}$ be a set of $(n+1)$-tuples of
convex functions on $(-a,a)$. 
Suppose that
there exist linear functions $\ell_\Phi$
such that the family
$$\{\Phi-\ell_\Phi\}=\left\{(\phi_0-\ell_\Phi,\ldots,\phi_n-\ell_\Phi)\right\}$$
is normal on $(-a,a)$.
Then one can take
$$\ell_\Phi(t)=\phi(0)+\phi'(0)t,\quad\mbox{where}\quad
\phi=\max_j\{\phi_0,\ldots,\phi_n\},$$
and $\phi'$ is the derivative from the right.
\end{lemma}

This is an immediate consequence from what was said before the lemma.

Now we return to our original setting: $f$ is a normal holomorphic map from
$\C^*$ to $\P^n$, and $F$ is a canonical representation of $f$
as in (\ref{1}), (\ref{co}).

We are going to state two
conditions for the curve $f$ to be normal.
We define
$$\N(t,F)=\max_{j=0}^n\N(t,g_j).$$
%
%
%
The derivative $
\N^\prime(t,F)$ is always understood as the derivative from the
right, so it takes only integer values, because $\N^\prime(t,g_j)$
takes only integer values.

Our first necessary condition of normality is a consequence of
Lemma~\ref{lemma3}:
\vspace{.1in}

{\em For every $a>0$ there exists $C(a)>0$ such that}
\begin{equation}\label{second}
\N(t,F)-\N(s,F)-\N'(s,F)(t-s)\leq C(a),\quad |t-s|<a.
\end{equation}
Condition (\ref{second}) is equivalent to
\begin{equation}\label{secondp}
\N(t,F)-\N(s,F)-\N'(s,F)(t-s)\leq C_1(1+(t-s)^2),
\end{equation}
for some $C_1>0$ and all $s,t$.

Our {\em second condition} is related to statement (ii)
of Lemma 1 and (\ref{cond2}):
\vspace{.1in}

{\em There exists $\delta>0$ and $C>0$ such that for every disc
(with respect to the metric (\ref{metric})) of radius $\delta$, centered
at a point $w\in\C^*$, the following
condition is satisfied: if the disc contains zeros of functions $g_j$ for 
$j\in I\subset\{0,\ldots,n\}$ then 
\begin{equation}\label{third}
N(|w|,F)\leq \max_{j\in I'}N(|w|,g_j)+C,
\end{equation}
where $I'$ is the complement of $I$ in $\{0,\ldots,n\}$.}
\vspace{.1in}

We will later prove that this condition is necessary for normality of $f$.

These two conditions of normality of a curve $f$ are formulated in
terms of parameters of formula (\ref{co}) for the homogeneous coordinates.
We give explicit expressions of the functions $N(r,g_j)$
in terms of these parameters:
$$N(r,g_j)=\log|A_j|+m_j\log r+\int_{[1,r]} n(1,t,g_j)\frac{dt}{t}.$$
Here $n(r_1,r_2,g_j)$
is the number of zeros of $g_j$ in the ring $r_1<|z|\leq r_2$.

\vspace{.1in}
Now we prove

\begin{theorem}
Conditions (\ref{second}) and (\ref{third})
are necessary and sufficient for normality of a curve $f$ with
homogeneous representation (\ref{1}), (\ref{co}).
\end{theorem}

Necessity of condition (\ref{second}) has already been proved in lemmas~\ref{lemma2}
and \ref{lemma3}.
To prove the rest we use the following

\begin{lemma}\label{lemma4}
Let $(f_k)$ be a sequence of holomorphic curves, $f_k(z)=f(\lambda_k z)$,
where $f$ has a homogeneous representation (\ref{1}), (\ref{co}).
Suppose that condition (\ref{second}) with $s=0$ holds for these curves
uniformly with respect to $k$, and 
set $$F_k(z)=(g_0(\lambda_kz),\ldots,g_n(\lambda_kz))=(g_{k,0},\ldots,g_{k,n}),$$
and 
\begin{equation}\label{mult}
h_k(z)=\exp(-\N(0,F_k))z^{-\N^\prime(0,F_k)}.
\end{equation}
Choose a subsequence on which $\N(t,h_kg_{k,j})$ tend to limits as $k\to\infty$,
and let $I$ be the set of indices $j$ for which the limit is finite,
and $I^\prime$ is the rest of the indices.

Then:

\noindent
a) $h_kg_{j,k}$ tend to limits, not identically equal to zero, for $j\in I$, and

\noindent
b) $h_kg_{j,k}$ tend to zero for $j\in I^\prime$. 
\end{lemma}

{\em Proof of the lemma.}
Condition (\ref{secondp}) implies that
\begin{equation}\label{upper}
\N(t,F_k)\leq C_1(1+t^2),
\end{equation}
for all $k$, $t$ and some $C_1>0$.

We first prove a). 
We define functions $h_k$ by (\ref{mult}). These functions are
holomorphic and zero-free because $\N^\prime(t,F)$ are integers.
For $j\in I$, the functions
$$\N(t,h_kg_{k,j})=\N(t,g_{k,j})-\N(0,f_k)-\N'(0,f_k)t$$
are convex, uniformly bounded on any interval,
the jumps of their derivatives are integers,
so the total jump of the derivatives is bounded on every interval.
Moreover, this total jump is at most a constant times
the length of the interval,
so we conclude that the $g_{k,j}$ have at most $C|b-a|$
zeros on every interval
$[a,b]$, so one can pass to the limit in formulas (\ref{co}), after multiplication of these
formulas by the $h_k$.
So the coordinates with $j\in I$ tend to non-zero limits, after choosing a subsequence.

Now we prove b), that is that remaining coordinates tend to zero.
We fix some $j\in I^\prime$ and will omit it
from the formulas, because all argument applies to any such coordinate.
We will also omit the index $k$ to simplify our formulas. So $hg=h_kg_{k,j}$.
We are going to prove that
$$B(x)=\log\max_{\theta}|(hg)(e^{x+i\theta})|\to -\infty,$$
uniformly for $|x|\leq \log 2.$
We assume for simplicity of formulas that $x=0$.
We represent $hg$ as a canonical product of
the form (\ref{co}), and denote by $n(r)$ the number of zeros on the
in the ring $\{ z:1\leq |z|\leq r\}$ if $r>1$ or
$\{ z:r<|z|<1\}$ if $r<1$, and $\n(t)=n(e^t).$
Then we have
$$B(0)\leq\log|A|+\int_1^\infty\log\left(1+\frac{1}{\xi}\right)dn(\xi)-
\int_0^1\log(1+\xi)dn(\xi).$$
Integrating by parts we obtain
$$B(0)\leq\log|A|+\int_1^\infty\frac{n(\xi)d\xi}{\xi(1+\xi)}+
\int_0^1\frac{n(\xi)d\xi}{(1+\xi)}.$$
Changing the variable $\xi=e^t$ gives
\begin{equation}\label{B}
B(0)\leq\log|A|+\int_{-\infty}^\infty\frac{\n(t)dt}{1+e^{|t|}}.
\end{equation}
Here $A$ is the number
from (\ref{co}) which depend on $gh=g_kh_k$.

On the other hand,
$$\N(x,hg)=\log|A|+mx+\int_{[0,x]}\n(t)dt.$$
By integrating (\ref{B}) by parts once more, we obtain
\begin{equation}\label{BB}
B(0)\leq \int_0^\infty(\N(t,hg)+\N(-t,hg))
\frac{e^t}{(1+e^t)^2}dt.
\end{equation}
Now $\N(0,hg)=\log|A|$, $\N$ is convex and $\N(t,hg)\leq C_1(1+t^2)$ in view
of (\ref{upper}). These conditions imply that
$$\N(t,hg)\leq 2(C_1+\sqrt{-C_1\log|A|})|t|+\log|A|.$$
Substituting this inequality to (\ref{BB}), we obtain
$B(0)\leq \log|A|+C_2\sqrt{-\log|A|}+C_3.$
As $A=A_k\to 0$, this completes the proof of the lemma.
\vspace{.1in}

Now necessity of condition (\ref{third}) follows immediately
because (\ref{third}) is now the same as (\ref{cond2}):
$N(|z_0|,h_kg_{k,j})=\log|h_kg_{k,j}(z_0)|+O(1)$ when $h_kg_{k,j}$ tends to a non-zero limit
as $k\to\infty$.

Sufficiency also immediately follows from Lemmas~\ref{lemma1} and \ref{lemma4}.
This completes the proof of the theorem.
\vspace{.1in}

Now we compare the result with Ostrowski's conditions.
His conditions are:
\vspace{.1in}

\noindent
a) The difference between the numbers of zeros and poles in any
ring $r_1<|z|<r_2$ is bounded uniformly with respect to $r_1,r_2$.
\vspace{.1in}

\noindent
b) The number of zeros of each coordinate in rings $r<|z|<2r$ is bounded,
\vspace{.1in}

\noindent
c) The distance between a zero of $g_0$ and a zero of $g_1$ is bounded
from below.
\vspace{.01in}

\noindent
d) There is a constant $C$ such that for each zero $w$ of $g_j$, 
we have $N(|w|,g_j)\leq N(|w|,g_{1-j})+C,\; j=0,1.$ 
\vspace{.1in}

It is easy to see that our condition (\ref{third}) with $n=1$ implies 
c) and d). When $|I|=2$, it is c) and when $|I|=1$, it is d).

Condition (\ref{second}) with $n=1$ implies b).
Condition a), which in our notation means
that $|\N^\prime(t,g_0)-\N^\prime(t,g_1)|$ is bounded,
is an easy consequence of (\ref{second}) and (\ref{third}),
when $n=1$.

However, for $n\geq 2$, conditions a) and b) do not have to hold.
Here is a simple example. Let 
$$g_0(z)=\prod_{j=0}^\infty(1-2^{-n}z),\quad g_1(z)=\prod_{j=0}^\infty
(1+2^{-n}z),$$
and take as $g_2$ any entire function with the property
$$T(r,g_2)=O(\log^{3/2}r),\quad r\to\infty.$$
It is easy to see that $\psi=f_0/f_1$ is a normal meromorphic function with
all limits $\lim_{\lambda\to\infty}\psi(\lambda z)$ non-constant, and
$$\log|g_3(z)|\leq o(\max\{\log|g_0(z)|,\log|g_1(z)|\}).$$
These properties imply that the curve $f$ with homogeneous coordinates
\newline
$(g_0,g_1,g_2)$ is normal, but $g_2$ can be chosen  so that
the number of its zeros in some rings $r<|z|<2r$ is unbounded. 

The author thanks Masaki Tsukamoto for his questions that stimulated
this work and Sergei Favorov for a useful discussion.

{\em Purdue University

West Lafayette, IN 47907

eremenko@math.purdue.edu}
\end{document}